\documentclass{article}
\usepackage[utf8]{inputenc}
\def\MR#1{\href{http://www.ams.org/mathscinet-getitem?mr=#1}{MR#1}}
\usepackage{amssymb}
\usepackage{hyperref}
\usepackage{amsmath}
\newtheorem{mylemma}{Lemma}
\newtheorem{myCo}{Corollary}
\newtheorem{mytheorem}{Theorem}
\title{The exact power law and Pascal pyramid
\footnote{This work was supported by the Russian Foundation for Basic Research, grant N~15-06-07402.}}
\author{Vladimir~V.~Bochkarev, Eduard~Yu.~Lerner}
\date{}

\begin{document}
\maketitle
\begin{abstract}       
Let $\omega_0, \omega_1,\ldots, \omega_n$ be a full set of outcomes (letters, symbols) and let positive  $p_i$, $i=0,\ldots,n$, be their probabilities ($\sum_{i=0}^n p_i=1$). Let us treat $\omega_0$ as a stop symbol; it can occur in sequences of symbols (we call them words) only once, at the very end. The probability of a word is defined as the product of probabilities of its letters. We consider the list of all possible words sorted in the non-increasing order of their probabilities. Let $p(r)$ be the probability of the $r$th word in this list. We prove that if at least one of ratios $\log p_i/\log p_j$, $i,j\in\{ 1,\ldots,n\}$, is irrational, then the limit $\lim_{r\to\infty} p(r)/r^{1/\gamma}$ exists and differs from zero; here $\gamma$ is the root of the equation $\sum_{i=1}^n p_i^\gamma=1$. Some weaker results were established earlier. We are first to write an explicit formula for this limit constant at the power function; it can be expressed (rather easily) in terms of the entropy of the distribution~$(p_1^\gamma,\ldots,p_n^\gamma)$.
\end{abstract}

{\it Key words:} Power law, Pascal pyramid, Entropy, Stirling formula, Gaussian approximation, Uniform distribution of sequences.

\section{Introduction. The statement of the main theorem}
The wide presence of power laws in real networks, in biology, economics, and linguistics can be explained in the framework of various mathematical models (see, for example, \cite{durett,mitzenmacher}). According to the Zipf law (\cite{Baayen}), in a list of wordforms ordered by the frequency of occurrence, the frequency of the $r$th wordform obeys a power function of~$r$ (the value $r$ is called the rank of the wordform). One can easily explain this law with the help of the so-called monkey model.

Assume that a monkey (see a detailed bibliographical review in~\cite{monkey}) types any of 26 Latin letters or the space on a keyboard with the same probability of $1/27$. We understand a word as a sequence of symbols typed by the monkey before the space. Let us sort the list of possible words with respect to probabilities of their occurrence (the empty word, whose probability equals $1/27$, will go first in this list followed by 26 one-letter words whose probabilities equal $1/27^2$ and then by $26^2$ possible two-letters words and so on). Evidently, the probability $p(r)$ of a word with the rank of $r$ satisfies the inequality
\begin{equation}
c_1 r^{-\alpha}< p(r)< c_2 r^{-\alpha},\quad \mbox{where $\alpha=\log 27/\log 26$, $c_1,c_2>0$.}
\label{eq:one}
\end{equation}
(here and below we use the symbol $\log$ if the base of the logarithm is not significant; but for the natural logarithm we use the symbol $\ln$).

Relatively recently inequality~(\ref{eq:one}) was generalized to the case of non\-equ\-ipro\-bable letters. Let $p_0$ be the probability that the monkey types the space, let $p_i$, $i=1,\ldots,n,$ denote probabilities of choosing the $i$th letter from the set of $n$ letters ($p_i>0$, $\sum_{i=0}^n p_i=1$), and let $p(r)$ be, as above, the probability of a word with a rank of $r$. Then, as is proved in \cite{Mit, we}, the following inequality analogous to~(\ref{eq:one}) takes place, namely, $\exists c_1,c_2:
0<c_1<c_2$ such that
\begin{equation}
c_1 r^{-\alpha}< p(r)< c_2 r^{-\alpha},\quad \mbox{where $\alpha=1/\gamma$}
\label{eq:two}
\end{equation}
and $\gamma$ is the root of the equation $\sum_{i=1}^n p_i^\gamma=1$ (evidently, $0<\gamma<1$). Note that inequality~(\ref{eq:two}) is equivalent to the boundedness of the difference $-\log p(r)-\alpha \log r$.

In the case, when the probability of each letter is not fixed, but depends on the previous one, words represent trajectories of a Markov chain with the absorbing state $\omega_0$ and transient states $\omega_1,\ldots,\omega_n$. Then the value $p(r)$ is the probability of the $r$th trajectory in the list of possible trajectories sorted in the non-increasing order of probabilities. In this case, the asymptotic behavior of $p(r)$ does not necessarily have a power order. Namely, in this case one of two alternatives takes place (\cite{canada,ourELA}). The first variant is that there exists the limit
\[\lim_{r\to\infty} \frac{-\log p(r)}{r^{1/m}}=c,\quad c>0,\]
where $m$ is some positive integer constant value that depends on the structure of the transition probability matrix and the structure of states, where the initial distribution of the Markov chain is concentrated. The second variant is that independently of the initial distribution there exists the following nonzero limit (the so-called \textit{weak} power law):
\[\lim_{r\to\infty} \frac{-\log p(r)}{\log r}.\]
This limit equals $1/\gamma$, where $\gamma$ is now defined with the help of the substochastic matrix~$P$ of transition probabilities where the row and the column that correspond to the absorbing state~$\omega_0$ are deleted. Namely, raising all elements of the mentioned matrix to the power of $\gamma$ would equate its spectral radius to 1.

These results were obtained independently in \cite{arXiv1,canada} and later refined in \cite{ourELA}. Namely, as appeared, the first alternative means the subexponential order of the asymptotics, i.e., in this case $\exists c_1,c_2: 0<c_1<c_2$ such that
\[
c_1 \exp(-c r^{-1/m})< p(r)< c_2 \exp(-c r^{-1/m}).
\]

The case of the second alternative is much more difficult. If the matrix~$P$ does not have the block-diagonal structure with coinciding powers such that raising elements of blocks to these powers makes the spectral radius equal 1, then one can replace the weak power law with a \textit{strong} one. Namely, in this case the asymptotic behavior of $p(r)$ has the power order, i.e. inequality~(\ref{eq:two}) is valid (with ``matrix'' $\gamma$ defined above). Therefore, inequality~(\ref{eq:two}) takes place in a ``typical'' case of letter probabilities.

However one more natural question still remains without an answer.

Inequality~(\ref{eq:two}) means that the asymptotic form has a power \textit{order} but does not imply the \textit{exact} power asymptotics. In a general case, as follows from the first example given in this section, useful properties can be established neither when letters in words are Markov-dependent no when they are independent. However, as we prove later in this paper, in a ``typical'' case, for words composed of independent letters, the asymptotic behavior of the function~$p(r)$ is exact power. The following theorem is valid.
\begin{mytheorem}[main]
\label{th:main}
Let at least one of ratios $\log p_i/\log p_j$, $i,j\in\{ 1,\ldots,n\}$, be irrational and let $\gamma$ be the root of the equation $\sum_{i=1}^n p_i^\gamma=1$. Then the limit 
\begin{equation}
\label{eq:thone}
\lim_{r\to\infty} \frac{p(r)^{-\gamma}/p_0^{-\gamma}}{r}
\end{equation} 
exists and equals $H(\mathbf{p}^\gamma)$, where $H(\mathbf{p}^\gamma)$ is the entropy of $\mathbf{p}^\gamma=(p_1^\gamma,\ldots,p_n^\gamma)$, i.e., 
$
H(\mathbf{p}^\gamma)=-\gamma \sum_{i=1}^n p_i^\gamma \ln p_i.
$
\end{mytheorem}

Here and below we always write the function under consideration in the numerator and do the norming (defined analytically) function in the denominator of the fraction, whose limit is to be calculated. In intermediate calculations it may be more convenient to do the opposite, but since this results only in the trivial raising of the limit constant to the power of $-1$, we sacrifice the convenience of calculations for the clarity of statements of results. Evidently, the theorem asserts that under certain assumptions there exists the nonzero limit $p(r)/r^{-\alpha}$ (where $\alpha=1/\gamma$) as $r\to\infty$.
It is equal to
$
p_0 H(\mathbf{p}^\gamma)^{-1/\gamma}.
$

Let us describe the structure of the rest part of the paper. In Section~2 we state the main theorem in terms of multinomial coefficients (of the Pascal pyramid). The proof of the theorem is reduced to the estimation of the limit behavior of the sum of these coefficients over some simplex. In Section~3 we prove an analog of this theorem with an integral in place of the sum. In this section we essentially use the Stirling formula which allows us to reduce calculations to the evaluation of a multivariate Gaussian integral. We establish an explicit formula for the determinant of the matrix of the quadratic form that defines the integrand. Finally, in Section~4 we prove that the ratio of the integral to the sum tends to 1. Here we use the general properties of the Riemann integral and uniformly distributed sequences. In conclusion we discuss possible generalizations and unsolved problems.

\section{Equivalent statements of the main theorem and the Pascal pyramid}

Let us first note that if $p_0\to 0$, then $\gamma\to 1$.
Reducing the nominator of fraction~(\ref{eq:thone}) by $p_0^{-\gamma}$, we write the following statement in this case:
\begin{mytheorem}[the case of $\gamma=1$]
\label{th:main2}
Let $p_i>0$ be the probability of the symbol~$\omega_i$, $i=1,\ldots,n,$ while $\sum_{i=1}^n p_i=1$ (there is no stop symbol). Assume that at least one of ratios $\log p_i/\log p_j$, $i,j\in\{ 1,\ldots,n\}$, is irrational. Let us consider all possible finite words (including the empty one), sort them in the non-increasing order of probabilities (we equate the probability of the empty word to 1 and calculate the probability of any other word as the product of probabilities of its letters). Let $p(r)$ be the probability of the $r$th word in the list (the word with the \textit{rank} of~$r$). Then the limit $\lim_{r\to\infty} p(r)/r^{-1}$ exists and equals $H^{-1}(\mathbf{p})$, where $H(\mathbf{p})$ is the entropy of the vector $\mathbf{p}=(p_1,\ldots,p_n),$ i.e., $H(\mathbf{p})=-\sum_{i=1}^n p_i \ln p_i$.
\end{mytheorem}

In the statement of Theorem~\ref{th:main2}, as well as in Theorem~\ref{th:main}, we use the bold font for the vector whose components are denoted by the same letter with the index ranging from~1 to~$n$. In what follows we use the bold font for analogous denotations without mentioning this fact.

One can easily see that Theorem~\ref{th:main2} is not just a particular case of Theorem~\ref{th:main}, but these theorems are equivalent. Namely, the replacement of probabilities $p_i^\gamma$ with new ones $p_i$ turns the general case into the particular one. Therefore, in what follows we neglect $p_0$, assuming (without loss of generality) that ${\sum_{i=1}^n p_i\!=\!1}$.

Fix some probability $q\in (0, 1]$ and denote by $Q(q)$ the rank of the last word whose probability is not less than $q$ in the list of all words sorted in the non-increasing order of their probabilities. Let us redefine the function $p(r)$ for noninteger $r$ as $p(r)=p(\lfloor r \rfloor)$ (here $\lfloor \cdot \rfloor$ is the integer part of a number). Evidently, functions $q=p(r)$ and $r=Q(q)$ ($q\in (0, 1]$, $r\ge 1$) are inverse (more exactly, quasi-inverse), 
namely, the graph of one of hyperbola-shaped, decreasing stepwise functions turns into another one when axes $r$ and $q$ switch roles (in the first case, $q$ is the argument and $r$ is the value, and vice versa in the second case).

It can be clearly seen
that $\lim_{r\to\infty} c\, p(r)/r^{-1}=1$ is equivalent to
\[\lim_{q\to 0} c^{-1} Q(q)/q^{-1}=1.\]
Therefore the equality in the assertion of Theorem~\ref{th:main2} is equivalent to that
\[
\lim_{q\to 0} Q(q)/q^{-1}=  H^{-1}(\mathbf{p}).
\]

Denote the logarithm of the denominator in the last fraction by $z=-\ln q$ (i.e., $q=e^{-z}$) and let $\tilde Q(z)=Q(e^{-z})$. In view of considerations in the above paragraph the equality in the assertion of Theorem~\ref{th:main2} is equivalent to that
\begin{equation}
\label{eq:mainlimit}
\lim_{z\to \infty} (\ln \tilde Q(z)- z)=-\ln H(\mathbf{p}).
\end{equation}

Recall the proof of inequality~(\ref{eq:two}) in~\cite{we}. It is reduced to the proof of the boundedness of the difference $\ln \tilde Q(z)- z$ for the introduced function $\tilde Q(z)$ with $z\ge 0$. Nonnegative values of $z$ form the definition domain of the function $\tilde Q(z)$ because $q\le 1 \Leftrightarrow z\ge 0$. For convenience we redefine the function $\tilde Q(z)$ by putting $\tilde Q(z)=0$ for $z<0$.

Let $a_i=-\ln p_i$. Considering all possible variants of the last letters in words, whose quantity equals the value of the function $\tilde Q$, we obtain the functional equation $\tilde{Q}(z)=\tilde{Q}(z-a_1)+\ldots+\tilde{Q}(z-a_n)+\chi(z)$, where $\chi$ is the Heaviside step (i.e., the function that vanishes with negative values of the argument and equals 1 with nonnegative values). For $z \geq M = \max\{a_1,\ldots,a_n\}$ we get the following recurrent correlation:
\begin{equation}
\label{eq:rec}
Q_n(z)=Q_n(z-a_1)+\ldots+Q_n(z-a_n),
\end{equation}
where $Q_n(z)=\tilde{Q}(z)+1/(n-1)$.

The equality $\sum_{i=1}^n p_i=1$ implies that the function $\mbox{const} \exp z$ satisfies Eq.~(\ref{eq:rec}). Since
the function $Q_n(z)$ takes a finite number of positive values
within the $[0,M]$ interval, there exist positive $c_1$ and $c_2$ such that
\begin{equation}
\label{eq:eneq}
c_1\exp z<Q_n(z)<c_2\exp z
\end{equation}
for all $0\le z\le M$.

Replacing terms in the right-hand side of the recurrent correlation~(\ref{eq:rec}) with their lower (upper) bounds, we extend the solution set of inequality~(\ref{eq:eneq}) to the domain $0\le z\le M+m$, where $m=\min\{a_1,\ldots,a_n\}$. Repeating this procedure several times, in a finite number of steps we prove that the inequality is valid for any arbitrarily large~$z$. Performing the logarithmic transformation of the inequality, we conclude that $\ln Q_n(z)-z$ is bounded, and then so is the difference $\ln\tilde Q(z)-z$.

Let us return to Theorem~\ref{th:main2}. As was mentioned above,
Theorem~\ref{th:main2} asserts (under certain assumptions) not only the  boundedness of $\ln\tilde Q(z)-z$ but also the validity of equality~(\ref{eq:mainlimit}).
Let us recall the combinatory sense of the function~$\tilde Q$; it is mentioned in~\cite{we}. Evidently, all words that contain $k_1$ letters of the 1st kind, $k_2$ letters of the 2nd kind, $\ldots$, and $k_n$ letters of the $n$th kind have one and the same probability of
$\mbox{Pr}(\mathbf{k}) = p_1^{k_1}\ldots p_n^{k_n}$ (i.e., $-\ln \mbox{Pr}(\mathbf{k})=\sum_{i=1}^n k_i a_i$), ranks of these words are consecutive. The quantity of such words is defined by the multinomial coefficient
\begin{equation}
\label{eq:M}
M(\mathbf{k}) = \frac{(k_1+\ldots+k_n)!}{k_1!\ldots k_n!}.
\end{equation}

Considering the nonnegative part of the $n$-dimensional integer grid and associating the point $(k_1\ldots,,k_n)$ with the number $M(k_1,...,k_n)$, we get one of variants of the Pascal pyramid. By the definition of the function $\tilde Q$ the value $\tilde Q(z)$ equals the sum of multinomial coefficients $M(\mathbf{k})$ over all integer vectors $\mathbf{k}$ that lie inside the $n$-dimensional simplex~$S(z)=\{\mathbf{x}: \mathbf{x}\ge 0,\ \sum_{i=1}^n a_i x_i\le z\}$:
\begin{equation}
\label{eq:tildeQ}
\tilde Q(z)=\sum_{\mathbf{k}\in S(z)} M(\mathbf{k}).
\end{equation}

As a result, we obtain one more equivalent statement of the main theorem, which we are going to prove.

\begin{mytheorem}[the multinomial statement]
\label{th:main3}
Let $a_i$, $i=1,\ldots,n,$ be arbitrary positive numbers such that at least one of ratios $a_i/a_j$, $i,j\in\{ 1,\ldots,n\}$ be irrational and $\sum_{i=1}^n p_i=1$, where $p_i=\exp(-a_i)$. Let a function $\tilde Q$ obey formula~(\ref{eq:tildeQ}). Then
\[
\lim_{z\to\infty} \frac{\tilde Q(z)}{\exp(z)}=H^{-1}(\mathbf{p}),
\]
where $H(\mathbf{p})=\sum_{i=1}^n a_i p_i$.
\end{mytheorem}
\section{The proof of an analog of Theorem~\ref{th:main3} with integration instead of summation}
\subsection{Reduction of the integration to the calculation of a Gaussian integral}
The function $M(k_1,\ldots,k_n)$ is defined for integer nonnegative vectors $\mathbf{k}$. Let us redefine it for noninteger vectors by replacing (in this case) $x!$ in Definition~(\ref{eq:M}) with $\Gamma(x+1)$. In what follows we use the denotation $M(x_1,\ldots,x_n)$ (or $M(\mathbf{x})$) for the corresponding function which is continuous for nonnegative $x_i$. Further we consider this function and study its properties only for such (nonnegative)~$x_i$.

In this section we prove the following theorem.
\begin{mytheorem}[on the integral]
\label{th:main4}
Let $a_i$, $i=1,\ldots,n,$ be arbitrary positive num\-bers such that $\sum_{i=1}^n p_i=1$, where $p_i=\exp(-a_i)$. Let a function $f(z)$ obey the formula~$f(z)=\int_{\mathbf{x}\,\in S(z)} M(\mathbf{x})\,d{\mathbf{x}}$, where $d{\mathbf{x}}=\prod_{i=1}^n dx_i$. Then
\[
\lim_{z\to\infty} \frac{f(z)}{\exp(z)}=H^{-1}(\mathbf{p}).
\]
\end{mytheorem}
{\bf Proof.}\ 
Let us first recall some evident properties of the integrand. Note that the existence of the (Riemann) integral of $f(z)$ over the compact set $S(z)$ evidently follows from the continuity of~$M(\mathbf{x})$ in the domain under consideration.

If all components of the vector $(x_1,\ldots,x_n)$, possibly, except one component~$x_i$, equal zero, then by definition we have $M(x_1,\ldots,x_n)\equiv 1$. Let us prove that otherwise the function $M(x_1,\ldots,x_n)$ is strictly increasing in~$x_i$. Since the gamma function is positive definite, it suffices to prove that in this case the partial derivative of $\ln M(x_1,\ldots,x_n)$ with respect to $x_i$ is positive. It equals
\[ (\ln \Gamma)'(x_1+\ldots+x_n+1)- (\ln \Gamma)'(x_i+1).\]
The positiveness of this difference follows from the fact that the function $(\ln \Gamma)'$ is increasing; this property, in turn, follows from the logarithmic convexity of the gamma function (it is well known~\cite{artin} that $(\ln \Gamma)''(x)=\sum_{i=0}^\infty \frac1{(i+x)^2}>0$ with~{$x>0$}).

The proved assertion implies that the function $M(\mathbf{x})$ attains its maximum in the domain $S(z)$ at the boun\-dary $\left\langle \mathbf{a},\mathbf{x} \right\rangle=z$, where $\left\langle \mathbf{a},\mathbf{x} \right\rangle=\sum_{i=1}^n a_i x_i$. Let us calculate the exact asymptotics of the maximal value of the function $M(\mathbf{x})$ in the domain $S(z)$ with $z\to\infty$. For the vector $\mathbf x$ we denote by $x$ the sum of its components and parameterize $\mathbf x$ by the value $x$ and ratios $q_i=x_i/x$:
\[x_i=q_i x, \ q_i\ge 0,\ i=1,\ldots,n,\quad \sum_{i=1}^n q_i=1.\]
Let us use one simplest corollary of the Stirling formula~\cite{artin}, namely, the fact that with a nonnegative argument the value of the difference $\ln \Gamma(x+1)-(x \ln(x)-x+\ln (x+1)/2)$ is bounded. We obtain that with any~$x>0$,
\begin{equation}
\label{eq:lnM}
\ln M(x_1,\ldots,x_n)=x H(\mathbf{q} )+ O(\ln(x+1)),\quad \mbox{where } H(\mathbf{q} )=-\sum_{i=1}^n q_i \ln q_i
\end{equation}
(this correlation is closely connected with the so-called entropy inequality for multinomial coefficients).

We seek for the maximum of this function with $z\to\infty$ under one additional condition (namely, the requirement that the maximum is attained at the boun\-dary) $\left\langle \mathbf{a},\mathbf{x} \right\rangle=z$, where $a_i=-\ln p_i$,  $0<p_i<1$, $\sum_{i=1}^n p_i=1$. Since $a_i>0$, we get $O(\ln(x+1))=O(\ln z)$. Moreover, the condition $\left\langle \mathbf{a},\mathbf{x} \right\rangle=z$ with mentioned $x_i$ gives the correlation
\begin{equation}
\label{eq:xz}
x=z H(\mathbf{q}; \mathbf{p})^{-1},\quad \mbox{where }H(\mathbf{q}; \mathbf{p})=\sum_{i=1}^n a_i q_i=-\sum_{i=1}^n q_i\ln p_i.
\end{equation}
Substituting this expression in~(\ref{eq:lnM}), we conclude that the maximum of $\ln M$ (accurate to $O(\ln z)$) is attained at a vector~$\mathbf{q}$ such that the fraction $H(\mathbf{q})/H(\mathbf{q};\mathbf{p}))$ takes on the maximal value. Recall that the difference $H(\mathbf{q};\mathbf{p})-H(\mathbf{q})$ takes on only nonnegative values and is called the Kullback--Leibler distance (divergence) $D(\mathbf{q}\,|\,\mathbf{p})$ between distributions $\mathbf{q}$ and $\mathbf{p}$~(see \cite{Kelbert}). The minimum of this difference is attained at only one value of $\mathbf{q}=\mathbf{p}$; evidently, an analogous assertion is also true for $H(\mathbf{q};\mathbf{p})/H(\mathbf{q})$: if $\mathbf{q}\neq \mathbf{p}$
\begin{equation}
\label{eq:less1}
H(\mathbf{q})/H(\mathbf{q};\mathbf{p})<1.
\end{equation}
Consequently, the maximum of the function $\ln M(\mathbf{x})$ in the domain $S(z)$ is attained (accurate to $O(\ln z)$) at the intersection of the hyperplane $\left\langle \mathbf{a},\mathbf{x} \right\rangle=z$ with the straight line $x_i=p_i x$, $i=1,\ldots n,$ where it equals $z+O(\ln z)$.

Let us now immediately prove Theorem~\ref{th:main4}. Note first that by using the L'Hopital rule we can reduce the proof to that of the formula obtained by differentiating the $f(z)/\exp(z)$ numerator and denominator with respect to~$z$ and to the proof of the equality
\begin{equation}
\label{eq:delta}
\lim_{z\to\infty} \frac{\hat f(z)}{\exp(z)}=H^{-1}(\mathbf{p}),
\end{equation}
where
$
\hat f(z)=\int_{\mathbf{x}\geq 0} M(\mathbf{x})\delta(z-\left\langle \mathbf{a},\mathbf{x} \right\rangle)\,d{\mathbf{x}},
$
and $\delta(\cdot)$ is the delta function.

Let~$\varepsilon$ be a real arbitrarily small positive value. Denote by $\Lambda_\varepsilon$ the sector consisting of points $\mathbf{x}$,
$x_i=q_i x$, $\sum_i q_i=1$, such that
\begin{equation}
\label{eq:eneqq}
p_i-\varepsilon<q_i<p_i+\varepsilon, \quad i=1,\ldots,n.
\end{equation}
With fixed $z$ on the hyperplane $\left\langle \mathbf{a},\mathbf{x} \right\rangle=z$ correlations~(\ref{eq:lnM}) and~(\ref{eq:xz}) take the form
\begin{equation}
\label{eq:lnm}
\ln M(\mathbf{x})=z H(\mathbf{q} )/H(\mathbf{q}; \mathbf{p})+ O(\ln(z)).
\end{equation}
Let us now strengthen inequality~(\ref{eq:less1}), namely, let us prove that if for $\mathbf{q}$ cor\-re\-la\-tions~(\ref{eq:eneqq}) are violated, then
\begin{equation}
\label{eq:HH}
H(\mathbf{q} )/H(\mathbf{q}; \mathbf{p})<1-C_1(\mathbf{p}) \varepsilon^2,
\end{equation}
where $C_1(\mathbf{p})$ is a positive constant independent of $\mathbf{q}$.

Since $H(\mathbf{q}; \mathbf{p})$ is a convex combination of $-\ln p_i$, it evidently is bounded:
\[0<\min\nolimits_{i} (-\ln p_i) \leq H(\mathbf{q}; \mathbf{p}) \leq \max\nolimits_{i} (-\ln p_i). \]
Consequently, formula~(\ref{eq:HH}) is equivalent to the inequality
\[
H(\mathbf{q}; \mathbf{p})-H(\mathbf{q})=D(\mathbf{q}\,|\,\mathbf{p})> C_2 (\mathbf{p}) \varepsilon^2.
\]
The latter correlation follows from the well-known property of the Kullback--Leibler divergence
\[
D(\mathbf{q}\,|\,\mathbf{p})\geq \frac14 \left( \sum_{i=1}^n |p_i-q_i| \right)^2
\]
(see, for example, lemma~3.6.10 in \cite{Kelbert}).

The proved inequality~(\ref{eq:HH}) (in view of formula~(\ref{eq:lnm})) implies that outside the domain $\Lambda_\varepsilon$ the function $M(\mathbf{x})$ is exponentially small in comparison to the maximal value inside the domain which equals $\exp(z)$. More precisely, with $x\not\in\Lambda_\varepsilon$, $\left\langle \mathbf{a},\mathbf{x} \right\rangle=z$ we get
\begin{equation}
\label{eq:eneqqq}
M(\mathbf{x})<\exp\left\{(1-C \varepsilon^2)z\right\}\quad \mbox{for some }~C>0.
\end{equation}
Note that the condition of the exponential smallness in comparison to $\exp z$ remains valid, even if $\varepsilon$ depends on $z$ and tends to 0 as~$z$ increases, though not too fast. In what follows we assume that
\[
\varepsilon=\varepsilon(z)=z^{-1/2+\delta}, \mbox{ where } \delta>0 \mbox{ is sufficiently small. }
\]

One can easily see that the same exponential upper bound as in~(\ref{eq:eneqqq}) also takes place
not only for the $M$ function but also for its integral over the domain which volume grows according to a
power law:
\[
\int_{\mathbf{x}\not\in \Lambda_{\varepsilon(z)},\mathbf{x}\geq 0} M(\mathbf{x})\delta(z-\left\langle \mathbf{a},\mathbf{x} \right\rangle)\,d{\mathbf{x}}\,<\,\exp\left\{(1-C \varepsilon^2)z\right\}
\]
with $z\to\infty$. Therefore in limit~(\ref{eq:delta}) we can treat $\hat f(z)$ as the integral
\begin{equation}
\label{eq:hatf}
\int_{\mathbf{x}\in \Lambda_{\varepsilon(z)}} M(\mathbf{x})\delta(z-\left\langle \mathbf{a},\mathbf{x} \right\rangle)\,d{\mathbf{x}}.
\end{equation}

Let us define the asymptotics~(\ref{eq:lnM}) of the function $M(\mathbf{x})$ in the domain~$\Lambda_{\varepsilon(z)}$ more precisely. Let us use the standard Stirling formula, namely, the fact that with $x\to\infty$ it holds $\ln \Gamma(x+1)=x \ln(x)-x+\ln(x)/2+\ln(2\pi)/2+R(x)$, where $0<R(x)<1/(12x)$. We obtain that in the domain~$\Lambda_{\varepsilon(z)}$,
\[
M(\mathbf{x})=\frac{1}{\sqrt{(2\pi)^{n-1}}}\exp\left\{ x H(\mathbf{q})+\ln(x)/2-\sum\nolimits_{i=1}^n \ln( x_i)/2+O(1/z)\right\}.
\]
Here, as usual, $x=\sum\nolimits_{i=1}^n x_i$; $q_i=x_i/x$. Therefore, we conclude that when considering the asymptotics of function~(\ref{eq:hatf}) we can treat $M(\mathbf{x})$ as follows:
\begin{equation}
\label{eq:widetildeM}
\widetilde M(\mathbf{x})=\frac{1}{\sqrt{(2\pi)^{n-1}}}\exp\left\{ x H(\mathbf{q})+\ln(x)/2-\sum\nolimits_{i=1}^n \ln( x_i)/2\right\}.
\end{equation}

In the latter formula we can write the exponent as
\begin{equation}
\label{eq:lefrig}
\left\{ \phantom{\sum\nolimits_{i=1}^n}\!\!\!\!\!\!\!\!\!\!\!\!\!\!\! \right\}=x\ln x+\ln(x)/2-
\sum\nolimits_{i=1}^n\left( x_i \ln x_i+\ln(x_i)/2\right).
\end{equation}
Let us write the Taylor expansion up to second-order terms near the maximum point in the plane $\left\langle \mathbf{a},\mathbf{x} \right\rangle=z$, i.e., near the point $\mathbf{x}'=\mathbf{p} z / H(\mathbf{p})$ (in what follows we denote by $x'_i$ coordinates of the point $\mathbf{x}'$ and do by $x'$ the sum of these coordinates which evidently equals $z H(\mathbf{p})^{-1}$).

First of all, note that
\[
\widetilde M(\mathbf{x}')=
\frac{1}{\sqrt{(2\pi x')^{n-1}\prod_{i=1}^n p_i}}\exp(z).
\]
One can easily calculate second derivatives of expression~(\ref{eq:lefrig})
\[
\frac{\partial^2}{\partial x_i\partial x_j}\left\{ \phantom{\sum\nolimits_{i=1}^n}\!\!\!\!\!\!\!\!\!\!\!\!\!\!\! \right\}=
\left\{
\begin{array}{ll}
x^{-1}-(2x^2)^{-1}, & \mbox{ if }i\neq j,\\
x^{-1}-(2x^2)^{-1}-x_i^{-1}+(2x_i^2)^{-1},& \mbox{ otherwise}
\end{array}
\right.
\]
(note that we do not use first derivatives in the Taylor expansion near the maximum point).

If $x\in\Lambda_\varepsilon$, then by formula~(\ref{eq:xz}) we have $x-x'=z (H(\mathbf{q}; \mathbf{p})^{-1}-H(\mathbf{p})^{-1})=z O(\varepsilon)$ (in the latter inequality we use the continuity of the function $H(\mathbf{q}; \mathbf{p})^{-1}$). Consequently,
\begin{equation}
\label{eq:ximxsi}
x_i-x'_i=x q_i-x' p_i=(x'+z O(\varepsilon))(p_i+O(\varepsilon))-x' p_i=z O(\varepsilon).
\end{equation}

In particular, with chosen $\varepsilon=\varepsilon(z)$ we have $|x_i-x'_i|=O(z^{1/2+\delta})$. We obtain that in the domain $\Lambda_{\varepsilon(z)}$,
\begin{eqnarray*}
\widetilde M(\mathbf{x})&=&\frac{1}{\sqrt{(2\pi x')^{n-1}\prod_{i=1}^n p_i}}\exp(z)\times \\
&&\exp\left\{\frac{\sum_{i,j=1}^n (x_i-x'_i)(x_j-x'_j)}{2 x'}-\frac{\sum_{i=1}^n(x_i-x'_i)^2}{2 p_i x'} +O(z^{-1/2+3\delta}) \right\}.
\end{eqnarray*}
Here the term $O(z^{-1/2+3\delta})$ contains both the remainder of terms of the series whose order exceeds 2 and the value of $O(z^{-1+2\delta})$ added by some omitted second-order terms. With $z\to\infty$ we can neglect the term of $O(z^{-1/2+3\delta})$. Therefore, in integral~(\ref{eq:hatf}) in place of $M(\mathbf{x})$ we should substitute the function $\widehat M(\mathbf{x})$ which differs from $\widetilde M(\mathbf{x})$ in the fact that its exponent does not contain the term of $O(z^{-1/2+3\delta})$.

Let us change variables  in the integral as follows: $y_i=(x_i-x'_i)/\sqrt{x'}.$ Since the degree of homogeneity of the delta-function equals~$-1$, we obtain that limit~(\ref{eq:delta}) coincides with
\begin{equation}
\label{eq:K}
\frac{1}{\sqrt{(2\pi)^{n-1}\prod_{i=1}^n p_i}}
\int_{\mathbb{R}^n} \delta(\left\langle \mathbf{a},\mathbf{y} \right\rangle)\exp\{-\left\langle \mathcal{B} \mathbf{y},\mathbf{y} \right\rangle/2\}\,d\mathbf{y},
\end{equation}
where $\mathcal{B}$ is the $n\times n$ matrix, all whose elements equal $-1$, except diagonal components which are greater by $1/p_i$.

\subsection{Calculation of the determinant}~

\begin{mylemma}
\label{lem:one1}
Let $n\ge 2$. Consider the $n\times n$ matrix $B$, where all nondiagonal elements equal 1, while $b_{ii}=1+k_i$. Then

1. The determinant of this matrix equals $\prod_{i=1}^n k_i \left(1+\sum_{j=1}^n 1/k_j\right)$.

2. The algebraic complement of the element with indices $(i,j), i\ne j,$ equals
\[-\prod_{\ell\in [n]\setminus\{i,j\}} k_\ell,\quad \mbox{where }[n]= \{1,\ldots,n\}.\]
\end{mylemma}

\begin{myCo}
The matrix $\mathcal{B}$ in formula~(\ref{eq:K}) is degenerate.
\end{myCo}
{\bf Proof of Lemma~\ref{lem:one1}.}\ 
Note that the first item of Lemma~\ref{lem:one1} defines the value of the algebraic complement of the diagonal element of such a matrix. Let us prove the theorem by induction.

With $n=2$ in the formula in item~2 we get the product over the empty set; it is accepted that this product equals 1. The formula in item~1 remains valid with $n=1$. In the induction step we assume that the formula in item~1 is proved for all dimensions less than~$n$ and has to be proved for the case when the dimension equals $n$, while the formula in item~2 is proved for all dimensions not greater than $n$ and has to be proved for the $(n+1)\times (n+1)$ matrix.

For proving item~1 we can use the expansion by the last row. Multiplying the algebraic complement by the diagonal element $k_n+1$, we get the sum
\begin{eqnarray}
&\prod_{i=1}^n k_i \left(1+\sum_{j=1}^{n-1} 1/k_j\right)
+ \prod_{i=1}^{n-1} k_i\left(1+ \sum_{j=1}^{n-1} 1/k_j \right)=& \nonumber
\\
\label{eq:sum1}
&\prod_{i=1}^n k_i \left(1+\sum_{j=1}^{n-1} 1/k_j\right)
+ \prod_{i=1}^{n-1} k_i+\sum_{j=1}^{n-1}\prod_{\ell\in [n]
\setminus\{n,j\}} k_\ell.&
\end{eqnarray}
The expansion by the entire last row, taking into account the induction hypothesis for item~2, make the third part in row~(\ref{eq:sum1}) vanish. First two terms in formula~(\ref{eq:sum1}) together give the desired sum.

In order to prove item~2, let us expand the determinant considered in this item  (algebraic complement of the element with $(i,j)$ indices of the
matrix~$B$ with the $(n+1)\times(n+1)$ dimension)
by the row whose number in the initial matrix of~$B$ was equal to~$j$. Generally speaking, for clarity, we use the same indices as in the numeration of the initial matrix. Since the algebraic complement considered in this item and the occurring algebraic complement for the element with indices $(j,i)$ (obtained by the expansion by a row of the determinant under consideration) have opposite signs, the value added by the element with indices $(j,i)$ equals
\[
-\prod_{r\in [n+1]\setminus \{i,j\}} k_r \left(1+\sum_{\ell\in [n+1]\setminus\{i,j\}} 1/k_\ell\right)
\]
(here we have used the induction hypothesis for item~1). The difference from the desired formula consists in the last term which equals (taking into account the first multiplier)
\[
-\sum_{\ell\in [n+1]\setminus\{i,j\}}\quad \prod_{r\in [n+1]\setminus\{i,j,\ell\}} k_r.
\]
It vanishes, when taking into account the contribution of the remaining $n-1$ elements in the $j$th row of the considered matrix.
\quad$\square$

\begin{mylemma}
\label{lem:twoo}
Let $B_1$ be the matrix mentioned in Lemma~\ref{lem:one1} (its dimension is $n\times n$, $n\geq 2$). Assume that $b_{ii}=1-1/p_i$, $i=1,\ldots n,$ where $p_i$ are arbitrary nonzero numbers. Denote by $B_2$ a matrix of the same dimension in the form $a^T a$, where $a=(a_1,\ldots a_n)$ is an arbitrary numeric row and $T$ is the transposition sign. Let $s$ be an arbitrary real number. Then
\[
\det\left( s B_2-B_1\right) =
\frac{s
   \left(\left(\sum _{i=1}^n a_i p_i \right)^2-\left(\sum
   _{j=1}^n p_j-1\right) \sum _{i=1}^n a_i^2
   p_i\right)}{\prod _{\ell=1}^n p_\ell}+\det(-B_1).
\]
\end{mylemma}

\begin{myCo}
\label{co:conc2}
If a vector $\mathbf{p}=(p_1,\ldots,p_n)$ satisfies additional constraints $p_i>0$, $\sum_{i=1}^n p_i=1$ (i.e., $-B_1=\mathcal{B}$), while $a_i=-\ln p_i$, then
\[
\sqrt{s\ \det\left( s^{-1} B_2-B_1\right)} =
\frac{H(\mathbf{p})}{\sqrt{\prod\nolimits _{\ell=1}^n p_\ell}}.
\]
\end{myCo}
{\bf Proof of Lemma~\ref{lem:twoo}.}\ 
By the differentiation rule for determinants, the derivative of the determinant of an $n\times n$ matrix equals the sum of determinants of~$n$ matrices such that in the $i$th one all elements of the $i$th row are replaced with their derivatives. We obtain that $\frac{\partial^2 \mbox{\rm det\,}\left( s B_2-B_1\right)}{\partial s^2}$ is the sum of determinants of matrices each one of which contains either the zero row or two various rows of the matrix~$B_2$. Since $\mbox{rank}\, B_2=1$, we get $\frac{\partial^2 \det\left( s B_2-B_1\right)}{\partial s^2}=0$.

Thus, $\det\left( s B_2-B_1\right)$ is a linear function of~$s$, whose free term evidently equals $\det(-B_1)$. It is clear that for calculating the coefficient $\det\left( s B_2-B_1\right)$ at $s$ it suffices to summate products of each element of the matrix $B_2$ by the algebraic complement of the corresponding element of the matrix~$-B_1$. If an element has indices $(i,j)$, $i\ne j,$ then by item~2 of Lemma~\ref{lem:one1} this product equals $a_i a_j p_i p_j/\prod _{\ell=1}^n p_\ell$.

Let`s explain the positive sign in the last formula. We calculate
an algebraic complement of the $-B_1$ matrix element. The
matrix has the $n\times n$ dimension, therefore the found algebraic
complement differs from the algebraic complement of the
corresponding $B_1$ matrix element for $(-1)^{n-1}$ times. According to
item~2 of Lemma~\ref{lem:one1}, the algebraic complement of the
corresponding $B_1$ matrix element is a ``minus'' product of $n-2$
multipliers~$k_i$. In the given case each of~$k_i$ factors is negative
(equals $-1/p_i$) which results in positive sign of the last formula in
the above paragraph.

Assume that this formula is valid for all $(i,j)$. Then we get the sum
\begin{equation}
\label{eq:two2}
\sum_{i,j=1}^n a_i a_j p_i p_j/\prod _{\ell=1}^n p_\ell=
\left(\sum _{i=1}^n a_i p_i \right)^2/\prod _{\ell=1}^n p_\ell.
\end{equation}
However by item~1 of Lemma~\ref{lem:one1} the algebraic complement of the diagonal element $b_{ii}$ of the matrix $-B_1$ equals
\begin{equation}
\label{eq:thr}
(-1)^{n-1}\left( \prod _{j:j\ne i} (-1/p_j)+\sum_{j:j\ne i}\prod _{\ell\not\in\{i,j\}}(-1/p_\ell)
\right)
\end{equation}
(here and below we omit the evident requirement that values of all indices belong to the set $[n]$).

Multiplying the first term in parentheses, i.e., $\prod _{j:j\ne i} (-1/p_j)$, by $(-1)^{n-1}a_i^2$ and summing over all~$i$, we get $\sum _{i=1}^n a_i^2 p_i/\prod _{\ell=1}^n p_\ell$. Let us multiply the resting term in parentheses~(\ref{eq:thr}) by $(-1)^{n-1}a_i^2$, sum over all~$i$, and subtract the value
\[\sum _{i=1}^n a_i^2 p_i^2/\prod _{\ell=1}^n p_\ell\]
from the obtained result (note that the subtrahend was ``illegally'' included in formula~(\ref{eq:two2})).
It gives the overall contribution of the second term in formula~(\ref{eq:thr}), which equals
\[
 -\sum_{j=1}^n p_j \sum _{i=1}^n a_i^2
   p_i /\prod _{\ell=1}^n p_\ell .
\]
Taking into account all the calculation elements of the
desired determinant allows completing the proof of Lemma~\ref{lem:twoo}.
\quad$\square$

For completing the proof of Theorem~\ref{th:main4} let us use Corollary~\ref{co:conc2}. Let us replace the variable in integral~(\ref{eq:K}) (as was proved earlier, this integral equals the limit considered in Theorem~\ref{th:main4})
$\delta(t)=\lim_{\sigma\to 0} \frac{1}{\sqrt{2\pi}\sigma}\exp\{-\frac{t^2}{2\sigma^2}\}$. Treating the limit multiplied by the coefficient at the exponent as a multiplier in the integral, we come to the limit of the Gaussian integral
\[
\lim_{\sigma\to 0}
\frac{1}{\sigma\sqrt{(2\pi)^{n}\prod_{i=1}^n p_i}}
\int_{\mathbb{R}^n} \exp\{-\left\langle (\sigma^{-2} B_2+\mathcal{B})\, \mathbf{y},\mathbf{y} \right\rangle/2\}\,d\mathbf{y},
\]
i.e.,
\[
\lim_{\sigma\to 0}
\frac{1}{\sigma \sqrt{\prod_{i=1}^n p_i\, \det\left( \sigma^{-2} B_2+\mathcal{B}\right)}}.
\]
Immediately applying Corollary~\ref{co:conc2}, we get desired~$H^{-1}(\mathbf{p})$.
\quad$\square$

\section{The ratio between the sum and the integral}
It remains to prove that under assumptions of Theorem~\ref{th:main3}, the ratio of the integral of the function~$M$ calculated over the domain $S(z)$ to the sum of values of this function at integer points of this domain tends to 1 as $z\to\infty$. For comparing the integral of the function and the sum of its values in the same domain one usually applies the Koksma--Hlawka inequality~(see \cite{niderraiter}). 
Note that usually one considers the integral over a fixed domain (as a rule, the cube $[0,1]^n$), whereas the domain in the case under consideration is varying. However, we intend only to prove the convergence of the fraction to 1 and do not need to estimate the asymptotic difference between the integral and the sum, which simplifies the task.

Evidently, it suffices to calculate the limit of the ratio for an arbitrary infinite increasing sequence $z_1,z_2,\ldots$ such that $z_i\to\infty$.

\begin{mytheorem}
\label{th:sumInt}
Let $\Omega_1,\Omega_2,\ldots$ be a sequence of Jordan measurable sets such that $\Omega_i\subset\Omega_{i+1}$ for all $i=1,2,\ldots$ Assume that $f(x)$, $x\in \Omega$, where $\Omega=\bigcup_i \Omega_i$, is an integrable and bounded on each of domains~$\Omega_i$ function such that $f(x)\ge 0$ and $\int f(x)\,d\, \Omega_i\to\infty$ as $i\to\infty$. Assume also that $K$ is a countable set of points from~$\Omega$ such that each of sets $K_i=K\cap \Omega_i$ is finite. Then if for any sufficiently small $\alpha>0$ there exists a partition of $\Omega$ onto a countable number of Jordan measurable sets $X_j=X_j(\alpha)$, $j=1,2,\ldots$, such that $\Omega_i=\bigcup_{j=1}^{n_i} X_j$ for some $n_i=n_i(\alpha)$, while
\begin{eqnarray}
\label{eq:supinf}
\sup_{x\in X_j}f(x)/\inf_{x\in X_j}f(x)<1+\alpha& &\mbox{starting with some } j,\\
\label{eq:limK}
\frac{|K\cap X_j|}{\mu X_j}\to 1& &\mbox{as } j\to\infty,
\end{eqnarray}
then in this case there exists the limit
\[
\lim_{i\to\infty} \frac{\int f(x)\,d\,\Omega_i}{\sum_{x\in K_i} f(x)}=1.
\]
\end{mytheorem}
{\bf Proof.}\ 
Evidently,
\[
\mu X_j \inf_{x\in X_j}f(x)\le \int f(x)\,d\,X_j \le  \mu X_j \sup_{x\in X_j}f(x).
\]
Therefore, in view of~(\ref{eq:supinf}) we conclude that starting with some $j$ it holds
\[
\frac{\int f(x)\,d\,X_j/\mu X_j}{\sum_{x\in K\cap X_j} f(x)/|K\cap X_j|}\in (1-\alpha,1+\alpha)
\]
with $\alpha<1$. In accordance with~(\ref{eq:limK}) we conclude that $\frac{\mu X_j}{|K\cap X_j|}\in (1-\alpha,1+\alpha)$ for all $j$, except a finite number of values of the index. Therefore, there exists $\ell$ such that with all~$j>n_\ell$,
\[
\frac{\int f(x)\,d\,X_j}{\sum_{x\in K\cap X_j} f(x)}\in ((1-\alpha)^2,(1+\alpha)^2).
\]
Representing this correlation as a double inequality and summing it over all~$j$ from $n_\ell+1$ to $n_i$, we obtain
\[
\frac{\int f(x)\,d(\Omega_i\setminus \Omega_\ell)}{\sum_{x\in K_i\setminus K_\ell} f(x)}\in ((1-\alpha)^2,(1+\alpha)^2)
\]
with $i>\ell$.

Note that by condition the numerator in the latter fraction (different from the integral $\int f(x)\,d\,\Omega_i$ by a constant value) tends to infinity. Then the same is true for the denominator. Note that the denominator differs from $\sum_{x\in K_i} f(x)$ by a constant value.

Therefore we conclude that all limit points of the sequence $\frac{\int f(x)\,d\,\Omega_i}{\sum_{x\in K_i} f(x)}$ lie inside the interval $((1-\alpha)^2,(1+\alpha)^2)$. Due to the arbitrariness of the choice of positive~$\alpha$ Theorem~\ref{th:sumInt} is proved.
\quad$\square$

\begin{myCo}[Completion of the proof of the main theorem]
\label{co:mainconc}
Let $f(z)$ be the function mentioned in assumptions of Theorem~\ref{th:main4}
and let $\tilde Q(z)$ obey fo\-r\-mu\-la~(\ref{eq:tildeQ}). Then if at least one of ratios $a_i/a_k$, $i,k\in\{ 1,\ldots,n\}$, $i\neq k,$ is irrational, then
\[
\lim_{z\to\infty} \frac{f(z)}{\tilde Q(z)}=1.
\]
\end{myCo}
{\bf Proof.}\ 
For clarity we denote by~$z$ the parameter that defines the boundary of the considered domain, and do by $\zeta$ the corresponding parameter of the hyperplane that contains a certain interior point $\mathbf{x}$ of this domain, i.e., $\zeta(\mathbf{x})=\left\langle \mathbf{a},\mathbf{x} \right\rangle$.

First of all, note that considerations in Section~3.1 imply that both in the sum and in the integral we can replace $S(z)$ with the domain
\[
\widehat \Lambda(z)= S(z)\cap \Lambda_{\varepsilon(\zeta)},\quad \mbox{ where }\varepsilon(\zeta)=\zeta^{-1/2+\delta},
\]
and replace the function $M(\mathbf{x})$ with that $\widetilde M(\mathbf{x})$ defined by formula~(\ref{eq:widetildeM}).
There\-fo\-re, we need to prove that
\begin{equation}
\label{eq:inttosumz}
\frac{\int_{\mathbf{x}\,\in \widehat \Lambda(z)}
\widetilde M(\mathbf{x})\,d{\mathbf{x}}}{\sum_{\mathbf{k}\,\in \widehat \Lambda(z)}\widetilde M(\mathbf{k})}\to 1
\end{equation}
(or that the difference of logarithms of the numerator and denominator tends to zero).

In view of Theorem~\ref{th:main4} the logarithm of the numerator in the latter fraction is a uniformly continuous function of $z$, while the logarithm of the denominator evidently is a nondecreasing function.
Therefore for proving the existence of the limit with $z\to\infty$ it suffices
to prove the existence of the limit for a sequence in the form $z_n=\kappa n$, $n=1,2,\ldots$,
where $\kappa$ is an arbitrarily small positive value
(as the difference between the numerator and denominator
of the logarithms in an arbitrary point slightly differs from the
value of difference in the nearest points $z_n$ in this sequence).
Namely, just for this fixed sequence we consider the ratio from the right-hand side of~(\ref{eq:inttosumz}).

In order to apply Theorem~\ref{th:sumInt}, for an arbitrary sufficiently small positive $\alpha$ we construct a partition of $\Lambda_{\varepsilon(\zeta)}$ onto domains $X_j$ satisfying assumptions of the theorem. Namely, we construct this partition by dividing of an infinite quantity of ``flapjacks'' located between neighboring hyperplanes in the form $\zeta({\mathbf{x}})=c_r$ and $\zeta({\mathbf{x}})=c_{r+1}$, $r=1,2,\ldots,$ where $c_{r+1}=c_r+\mbox{Const},$ onto a finite number of domains~$X_j$.

Evidently, for any $\alpha\le 2\kappa$ we can choose a sequence $c_r$ such that
\[
c_{r+1}-c_r=\mbox{Const}< \frac{\alpha}{2};
\quad \mbox{for any }n\
\exists\, r: z_n=c_r.
\]
To this end, it suffices to put $c_r=\mbox{Const}\,r$, where $\mbox{Const}= \kappa/ \lceil \frac{2\kappa}{\alpha} \rceil$
(here $\lceil \cdot \rceil$ is an upward rounding to the nearest integer).

Let $C_r=\{\mathbf{x}: c_r\le\zeta(\mathbf{x})< c_{r+1} \}$. Denote by $F_r$ the $r$th ``flapjack'' $C_r \cap \Lambda_{\varepsilon(\zeta)}$. We are going to ``cut'' $F_r$ onto a finite number of domains~$X_j$. We numerate the countable number of domains $X_j$, $j=1,2,\ldots,$ so as to make domains $X_j$ obtained by ``cutting'' $F_r$ with the least~$r$ have lesser numbers, while the order of numbering inside the partition of $F_r$ plays no role.

Since $\varepsilon(\zeta)=o(\zeta)$, with $\mathbf{x}, \mathbf{y}\in F_r$ it holds: $y_i=x_i+o(x_i)$ (cf. with~(\ref{eq:ximxsi})). Consequently, with $r\to\infty$ we get $\ln y_i-\ln x_i \to 0$ and $\ln y-\ln x \to 0$.

By formula~(\ref{eq:widetildeM}),
\[
\ln \widetilde M(\mathbf{x})=\mbox{const}+
g(\mathbf{x}) +\ln(x)/2-\sum\nolimits_{i=1}^n \ln( x_i)/2,
\]
where $g(\mathbf{x})=x H(\mathbf{q})=\sum_{i=1}^n x_i \ln x_i-\left(\sum_{i=1}^n x_i\right)\ln \left( \sum_{j=1}^n x_i\right).$
We get
\[
\mbox{grad} g=\ln \mathbf{q}=(\ln q_1,\ldots,\ln q_n),\quad
\frac{\partial^2 g}{\partial x_i\partial x_j}=O(1/x),\ i,j\in\{1,\ldots,n\}.
\]
Using expansion in a series $\ln \widetilde M$ with evaluation of the second
order terms and considerations of the previous paragraph we obtain the following important observation.
If $\mathbf{x}, \mathbf{y}\in F_r$~and
\begin{equation}
\label{eq:cr}
|x_i-y_i|=o(\sqrt{x})=o(\sqrt{c_r}), \quad i=1,\ldots,n,
\end{equation}
then with sufficiently large $r$ it holds
\[
\ln \widetilde M(\mathbf{x})-\ln \widetilde M(\mathbf{y})<
|\left\langle \ln \mathbf{q}\,,\mathbf{x}-\mathbf{y} \right\rangle|+\alpha/100.
\]

Since $\mathbf{q}\to \mathbf{p}$ as $r\to\infty$, with sufficiently large~$r$ it holds
\[
|\left\langle \ln \mathbf{q}\,,\mathbf{x}-\mathbf{y} \right\rangle|\,<\,
|\left\langle \mathbf{a}\,,\mathbf{x}-\mathbf{y} \right\rangle|+\alpha/200\,<\,
(c_{r+1}-c_r)+\alpha/200.
\]
As a result, we obtain that with sufficiently small $\alpha$, starting with some~$r$, it holds
\[
\widetilde M(\mathbf{x})/\widetilde M(\mathbf{y})<1+\alpha .
\]

Therefore, dividing $F_r$ onto domains $X_j$ so as to fulfill correlation~(\ref{eq:cr}) for all points $\mathbf{x},\mathbf{y}$ that belong to one domain, we guarantee the validity of assumption~(\ref{eq:supinf}) in Theorem~\ref{th:sumInt}. Note that it suffices to fulfill condition~(\ref{eq:cr}) for all indices~$i$ except one, because the validity of this condition for the rest index follows from the fact that $\mathbf{x},\mathbf{y}\in C_r$.

Finally, let us use the irrationality of $a_{i^*}/a_{k^*}$ for some $i^*\neq k^*$. Let us denote by $I_{k^*}$ the set $\{ 1,\ldots,n\}\setminus\{ k^* \}$ and do by $I_{i^*k^*}$ the set $\{ 1,\ldots,n\}\setminus\{i^*,k^* \}$. We are going to prove that defining domains $X_j$ by inequalities
\begin{equation}
\label{eq:IK}
l_{ji}\le x_i<L_{ji},\ i\in I_{k^*}, \qquad \mbox{where } L_{ji}-l_{ji}>\mbox{const } c_r^{1/2-\delta},
\end{equation}
we fulfill condition~(\ref{eq:limK}) (with $K=\mathbb{Z}^n$). Here, as usual, $\delta$ is a sufficiently small real positive value, though in this case we can choose $\delta$ as any number in the interval $(0,1/2)$
(roughly speaking, it`s sufficient that the radius of the pieces $X_j$
used to divide ``flapjacks'' $F_r$ tend to infinity at~$r\to\infty$).

Evidently, we can divide ``almost all'' $F_r$ onto domains~$X_j$ so as to simul\-ta\-neo\-usly fulfill inequalities~(\ref{eq:cr}) and conditions~(\ref{eq:IK}) on~$l$ and~$L$ (the remaining ``cuttings'' on the edges of the domain $F_r$ which occur due to the inconsistency between the inequality $l_{ji}\le x_i<L_{ji},\ i\in I_{k^*}$ and the definition of the boundary of the domain $\Lambda_{\varepsilon(z)}$ are asymptotically small).

Evidently, $\mu X_j=\prod_{i\in I_{k^*}} (L_{ji}-l_{ji})\times (c_{r+1}-c_r)/a_{k^*}$. Since the difference $(L_{ji}-l_{ji})$ grows as $j\to\infty$, \textit{the asymptotics} of the number of ways for choosing integer $x_i$ such that
$l_{ji}\le x_i<L_{ji}$ for $i\in I_{i^*k^*}$ coincides with $\prod_{i\in I_{i^*k^*}} (L_{ji}-l_{ji})$. Here and below we understand the asymptotics as a function of~$j$ such that the ratio of the considered quantity to this function tends to 1 as $j\to\infty$. In order to complete the proof of Corollary~\ref{co:mainconc}, it remains to prove the following lemma.

\begin{mylemma}
\label{lem:finlemma}
Let the ratio $a_{i*}/a_{k^*}$ be irrational and $(L_{ji^*}-l_{ji^*})\to\infty$. Assume also that the ratio $(c_{r+1}-c_r) /a_{k^*}$ equals a constant value lesser than 1 which is independent of~$r$. Then for fixed $x_i,\ i\in I_{i^*k^*},$ the asymptotics of the number of ways to choose integer $x_i$, $i\in \{i^*,k^*\},$ such that $l_{ji^*}\le x_{i^*}<L_{ji^*}$ and $c_r\le\zeta({\mathbf x})<c_{r+1}$ simultaneously, equals $(L_{ji}-l_{ji})\times (c_{r+1}-c_r)/a_{k^*}$.
\end{mylemma}
{\bf Proof of Lemma~\ref{lem:finlemma}.}\ 
In what follows we need standard denotations for the fractional part $\{\cdot\}$, floor $\lfloor\cdot\rfloor$, and ceil $\lceil\cdot\rceil$ of a number.

Let $c'=\sum_{i\in I_{i^* k^*}} a_i x_i$,
$d_r=(c_r-c')/a_{k^*}$,
$D_r=(c_{r+1}-c')/a_{k^*}$, and $\theta=a_{i*}/a_{k^*}$. The condition $c_r\le\zeta({\mathbf x})<c_{r+1}$ is equivalent to that
\begin{equation}
\label{eq:fincond}
\theta\, x_{i^*}+x_{k^*}\in [d_r,D_r).
\end{equation}

If the difference $D_r-d_r$ (it equals $(c_{r+1}-c_r) /a_{k^*}$) is less than~1
(this inequality is obviously holds for sufficiently small $\alpha$) then with fixed $x_{i^*}$ the integer value $x_{k^*}$ satisfying condition~(\ref{eq:fincond}) is defined uniquely, provided that it exists. Therefore, we need to estimate the quantity of values $x_{i^*}$ in the interval $[l_{ji^*},L_{ji^*})$ such that $\{ \theta\, x_{i^*}\}\in [\{d_r\},\{ D_r\})$; here the latter correlation is understood in the sense of an interval on the unit circle, and the length of the considered interval is independent of~$r$.

Recall the definition of a well-distributed sequence~\cite[section 1.5]{niderraiter}.

``Let $(y_n)$ $n=1,2,\ldots,$ be a sequence of real numbers. For integers $N\geq 1$ and $k\geq 0$ and a subset $E$ of $[0,1)$, let $A(E,N,k)$ be the number of terms among $\{y_{k+1}\},\{y_{k+2}\},\ldots,\linebreak
\{y_{k+N}\}$ that are lying in~$E$.

The sequence $(y_n)$ $n=1,2,\ldots,$ is said to be \textit{well-distributed mod 1} if for all pairs $a$, $b$ of real numbers with $0\le a<b\le 1$ we have
\[
\lim_{N\to\infty} A([a,b);N,k)/N=b-a\quad\mbox{uniformly in }k=0,1,2,\ldots.
\]

Example. The sequence $(n\theta)$ $n=1,2,\ldots$ with $\theta$ irrational is well-distributed {mod 1}.''

The latter fact would have proved Lemma~\ref{lem:finlemma}, if the interval of the unit circle $[\{ d_r\},\{D_r\})$ was independent of~$r$. Let us clarify this property in the case of the inequality $\{D_r\}>\{ d_r\}$. In what follows we always assume that this inequality is valid (evidently, as in the definition of the well-distribution property, this leads to no loss of generality). Really, if $k$ equals $\lceil l_{ji^*} \rceil-1$ and $N$ does the difference $\lfloor L_{ji^*} \rfloor-\lceil l_{ji^*} \rceil+1$, then we obtain the uniform in $j$ convergence
\begin{equation}
\label{eq:lastlimit}
\lim_{N\to\infty} \frac{A(\,[\{ d\},\{D\});N,\lceil\, l_{ji^*} \rceil-1)}{N}=D-d,
\end{equation}
which is equivalent to the assertion of the Lemma with the fixed of $\{ d_r\}$, $\{D_r\}$.

Note that if with fixed~$j$ equality~(\ref{eq:lastlimit}) is valid for any subinterval in $[0,1)$, then we say that the corresponding sequence is \textit{uniformly distributed modulo~1}. This property follows from the property of the \textit{well-distribution modulo~1}. It is well known that (see \cite{niderraiter}[section 2.1]) for any sequence uniformly distributed modulo~1 the convergence is uniform with respect to all subintervals in $[0,1)$. Consequently, we get the uniform in~$j$ convergence
\[
\lim_{N\to\infty} \frac{A(\,[\{ d_r\},\{D_r\});N,\lceil\, l_{ji^*} \rceil-1)}{N}=\Delta,
\]
where the \textit{constant} $\Delta$ equals $D_r-d_r$. Therefore,
\[
\lim_{j\to\infty} \frac{A(\,[\{ d_r\},\{D_r\});\lfloor L_{ji^*} \rfloor-\lceil l_{ji^*} \rceil+1,\lceil\, l_{ji^*} \rceil-1)}{(L_{ji^*}-l_{ji^*})\,\Delta}=1,
\]
which coincides with the lemma assertion in a general case.
\quad$\square$

\noindent
$\square$

\section{Conclusion}
We have proved that in the monkey model the probability of words in the sorted list has the \textit{exact} power asymptotics, provided that the ratio of logarithms of probabilities of certain letters is irrational.

Note that this condition is not only sufficient, but also necessary. Really, otherwise logarithms of probabilities $a_i=-\ln p_i$, $i=1,\ldots,n,$ allow the repre\-sen\-ta\-tion $a_i=m_i v$, where $m_i$ are natural numbers and $v$ is independent of~$i$. In this case formula~(\ref{eq:rec}) defines a linear recurrent correlation on a grid with the step of $v$. This does not affect the initial constancy of the function~$\tilde Q$ in cells of the grid with the mentioned step with any value of the argument.

It should be noted that using the expression for terms of
linear recurring sequences via the corresponding powers of
roots of the characteristic equation allows clear analysis of rate
of convergence to the power law of the function~$\tilde Q$ (with a step
of~$v$ on the grid) in this degenerate case. It would be more
interesting to conduct such studies for more general case to
which the main theorem of this paper is devoted.

A generalization of results obtained in this paper to the case of the Markov dependence is of even more interest. In this case an analog of the vector $\mathbf{p}^\gamma$ is a substochastic matrix of transition probabilities where the row and column that correspond to the absorbing state are deleted, and all elements of this matrix are raised to a power of~$\gamma$ such that its spectral radius equals 1. Denote this matrix by $\mathbf{P}^\gamma$. In the case considered above all rows of the matrix $\mathbf{P}^\gamma$ coincide with $\mathbf{p}^\gamma$. In a typical case, when the strong power law takes place (see Introduction), the matrix $\mathbf{P}$ is irreducible and the matrix transposed with respect to $\mathbf{P}^\gamma$ has a positive eigenvector that corresponds to the unit eigenvalue. Let us norm this vector so as to make the sum of its components equal 1 and denote the result by~$\mathbf{w}$. In the case of the Markov chain with the transition probability matrix~$\mathbf{P}^\gamma$ this vector defines an ergodic distribution.

If all rows of the considered matrix coincide with $\mathbf{p}^\gamma$, then one can easily see that $\mathbf{w}$ coincides with $\mathbf{p}^\gamma$. It is possible that in the case of the Markov dependence with the irreducible matrix $\mathbf{P}$, an analog of Theorem~\ref{th:main} takes place, where the role of the vector $\mathbf{p}^\gamma$ is played by $\mathbf{w}$. We are going to verify this conjecture in our future research.


\begin{thebibliography}{99}
\bibitem{artin}
{Artin, E.}: 
\newblock \emph{The gamma function}.
\newblock Translated by Michael Butler. Athena Series: Selected Topics in
  Mathematics. Holt, Rinehart and Winston, New York-Toronto-London, 1964.  
\MR{0165148}
\bibitem{Baayen}
{Baayen, R.~H.}:
\newblock \emph{Word frequency distributions}. Text, Speech and Language
  Technology, Vol.~\textbf{18}.
\newblock Kluwer Academic Publishers, Dordrecht.
\newblock With 1 CD-ROM (Windows and Unix), 2001.
\MR{1855236}
\bibitem{we}
{Bochkarev, V.~V.} {and} {Lerner, E.~Yu.}:
\newblock The {Z}ipf law for random texts with unequal letter probabilities and
  the {P}ascal pyramid.
\newblock \textit{ Izv. Vyssh. Uchebn. Zaved. Mat.\/}~\textbf{56},  (2012), no.~12, 30--33.
\MR{3137107}
\bibitem{arXiv1}
{Bochkarev, V.~V.} {and} {Lerner, E.~Yu.}: 
Zipf and non-Zipf laws for homogeneous Markov chain.  Preprint, 2012. 
\href{http://arxiv.org/abs/1207.1872}{arXiv:{1207.1872}}.
\bibitem{ourELA}
{Bochkarev, V.~V.} {and} {Lerner, E.~Yu.}:
\newblock Strong power and subexponential laws for an ordered list of
  trajectories of a {M}arkov chain.
\newblock \emph{Electron. J. Linear Algebra\/}~\textbf{27},  (2014), 534--556.
\MR{3266165}
\bibitem{Mit}
{Conrad, B.} {and} {Mitzenmacher, M.}:
\newblock Power laws for monkeys typing randomly: the case of unequal
  probabilities.
\newblock \emph{IEEE Trans. Inform. Theory\/}~\textbf{50}, (2004), no.~7, 1403--1414.
\MR{2095846}
\bibitem{durett}
{Durrett, R.}:
\newblock \emph{Random graph dynamics}.
\newblock Cambridge Series in Statistical and Probabilistic Mathematics.
  Cambridge University Press, Cambridge, 2007.
\MR{2271734}
\bibitem{canada}
{Edwards, R.}, {Foxall, E.}, {and} {Perkins, T.~J.}: 
\newblock Scaling properties of paths on graphs.
\newblock \textit{ Electron. J. Linear Algebra\/}~\textbf{23}, (2012), 966--988.
\MR{3007200}
\bibitem{niderraiter}
{Kuipers, L.} {and} {Niederreiter, H.}:
\newblock \textit{Uniform distribution of sequences}.
\newblock Wiley-Interscience [John Wiley \& Sons], New York-London-Sydney.
\newblock Pure and Applied Mathematics, 1974. 
\MR{0419394}
\bibitem{mitzenmacher}
{Mitzenmacher, M.}:
\newblock A brief history of generative models for power law and lognormal
  distributions.
\newblock \textit{Internet Math.\/}~\textbf{1}, (2004), no.~2, 226--251.
\MR{2077227}
\bibitem{monkey}
{Perline, Richard} {and} {Perine, Ron}:
\newblock Two Universality Properties Associated with the Monkey Model of Zipf's Law. 
\textit{ Entropy\/}~\textbf{18},  (2016), no.~3, 89; \href{http://www.mdpi.com/1099-4300/18/3/89}{doi:10.3390/e18030089}.
\bibitem{Kelbert}
{Suhov, Yu.} {and} {Kelbert, M.}: 
\newblock \textit{Probability and statistics by example. {II}}.
\newblock Cambridge University Press, Cambridge.
\newblock Markov chains: a primer in random processes and their applications, 2008.
\MR{2423218}
\end{thebibliography}
\end{document}